\documentclass[11pt]{article}%
\usepackage{amsfonts}
\usepackage{amsmath}
\usepackage{amssymb}
\usepackage{graphicx}%
\providecommand{\U}[1]{\protect\rule{.1in}{.1in}}
\newtheorem{theorem}{Theorem}

\newtheorem{corollary}[theorem]{Corollary}

\newenvironment{proof}[1][Proof]{\noindent\textit{#1.} }{}
\begin{document}

\begin{center}
{\Large Constructing a quadrilateral inside another one }
\end{center}

\vspace*{0.2 in}

\begin{flushright}
J. Marshall Ash
\\
DePaul University
\\
Chicago, IL 60614 \vspace*{0.25in}

Michael A. Ash
\\
Department of Economics
\\
University of Massachusetts Amherst
\\
Amherst, MA 01003 \vspace*{0.25in}

Peter F. Ash
\\
Cambridge College
\\
Cambridge, MA 02138\vspace*{0.25in}
\end{flushright}

\section{\label{s:1}The quadrilateral ratio problem}

The description of Project 54 in 101 Project Ideas for the Geometer's
Sketchpad \cite{Key} reads (in part):

\qquad\parbox[1in]{4.5in}{On the Units panel of Preferences, set
Scalar Precision to hundredths. Construct a generic quadrilateral
and the midpoints of the sides. Connect each vertex to the
midpoint of an opposite side in consecutive order to form an inner
quadrilateral. Measure the areas of the inner and original
quadrilateral and calculate the ratio of these areas. What
conjecture are you tempted to make? Change Scalar Precision to
thousandths and drag until you find counterexamples.}

A figure similar to the following accompanies the project.
\begin{center}
\includegraphics[
height=2.911in,
width=4.0318in
]%
{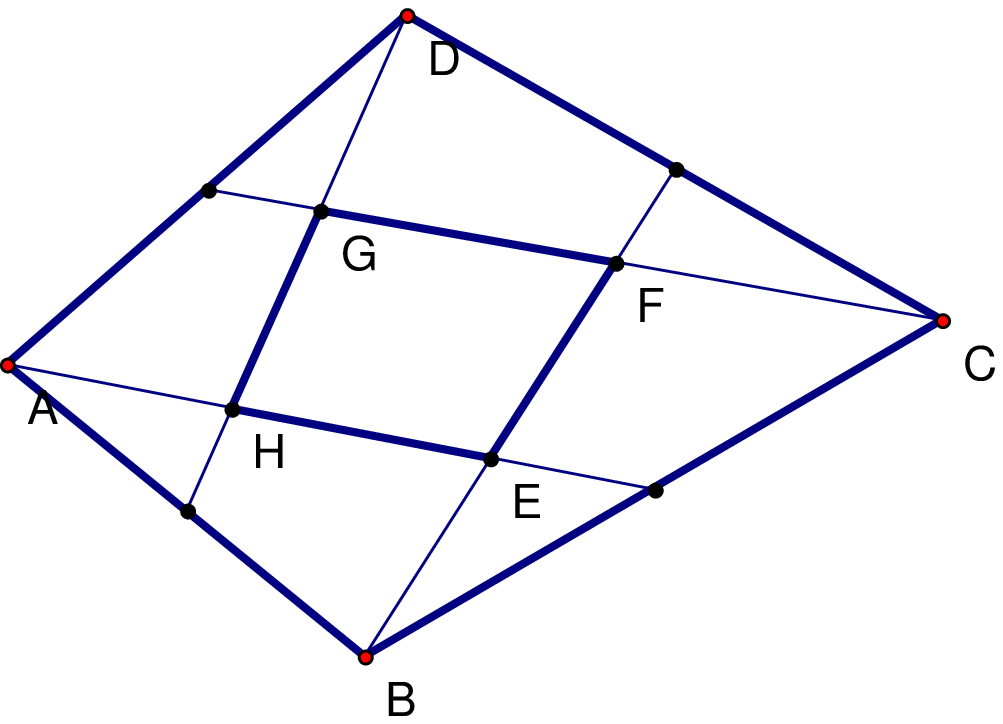}%
\end{center}
Let $r$ be the ratio of the quadrilateral areas,%
\begin{equation}
r=\frac{\text{area}\left(  EFGH\right)  }{\text{area}\left(  ABCD\right)
}\text{.} \label{ratio}%
\end{equation}
The tempting conjecture is that $r=1/5$. In Theorem 1 we show that this is
true in case the original quadrilateral is a parallelogram. However, the
conjecture is false in general. Instead, the ratio can be any real number in
the interval $(1/6,1/5]$. This is our Corollary \ref{c:2}. (Since
$1/6<0.17<1/5$, it is possible to find counterexamples even with the Sketchpad
Scalar Precision set to hundredths; however, it takes very industrious
dragging to find them.)

Suppose now that instead of a quadrilateral we had a triangle. Of course,
joining each vertex to the opposite midpoint would not yield an inner
triangle, since the three lines are medians, which are concurrent in a point.
To look for an analogous result for a triangle, we can look for points which
are not midpoints, but rather divide each side a ratio $\rho$ of the distance
from one point to the next, $0<\rho<1$. For definiteness, we assume that
\textquotedblleft next point\textquotedblright\ in this definition is based on
movement in a counterclockwise direction. We call these points $\rho
$\textit{-points}. It turns out that for a given $\rho$, the ratio of the area
of the inner triangle to the area of the outer triangle is constant
independent of the initial triangle and is given by $\frac{\left(
2\rho-1\right)  ^{2}}{\rho^{2}-\rho+1}$. Note that when $\rho=1/2$, this
reduces to $0$, which provides a convoluted proof that the medians of a
triangle are concurrent. When $\rho=1/3$, the area ratio is $1/7$, which the
Nobel Prize-winning physicist Richard Feynman once proved, though he was
probably not the first to do so. The result for general $\rho$ is known, and a
proof is given in \cite{DeV}, along with the Feynman story.

Inspired by this result, we will study the quadrilateral question for $\rho
$-points. Let $ABCD$ be a convex quadrilateral and let $N_{1},N_{2},N_{3},$
and $N_{4}$ be chosen so that $N_{1}$ is the $\rho$-point of $BC$, $N_{2}$ is
the $\rho$-point of $CD$, $N_{3}$ is the $\rho$-point of $DA$, and $N_{4}$ is
the $\rho$-point of $AB$. For fixed $\rho$ $\left(  0<\rho<1\right)  $ connect
each vertex of $ABCD$ to the $\rho$-point of the next side. ($A$ to $N_{1}$,
$B$ to $N_{2}$, $C$ to $N_{3}$, and $D$ to $N_{4}$.) The intersections of the
four line segments form the vertices of a convex quadrilateral $EFGH$.%
\begin{center}
\includegraphics[
height=2.469in,
width=4.8231in
]%
{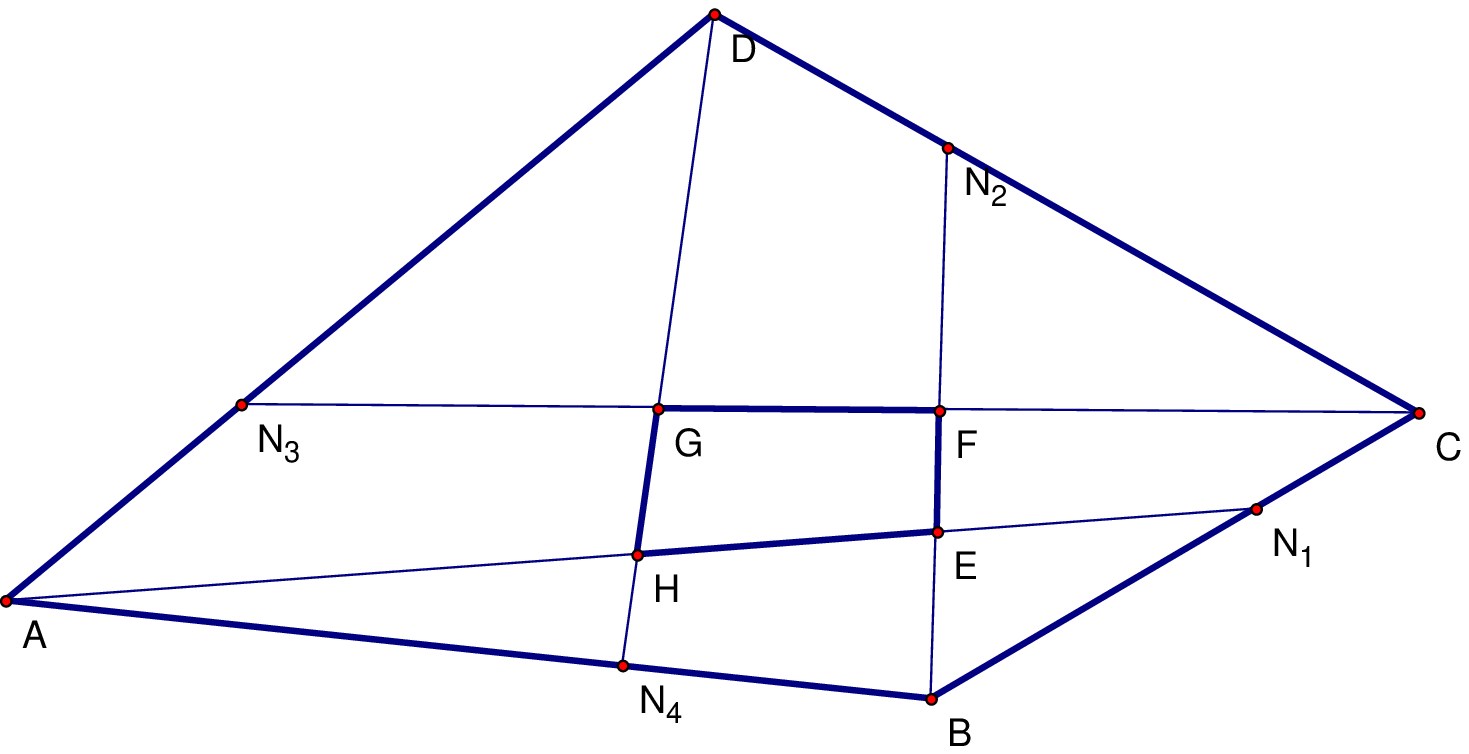}%
\end{center}

Define the area ratio
\[
r(\rho,ABCD)=\frac{Area\left(  EFGH\right)  }{Area\left(  ABCD\right)  }.
\]
Theorem \ref{t:2} below states that as $ABCD$ varies, the values of $r\left(
\rho,ABCD\right)  $ fill the interval%
\begin{equation}
\left(  m,M\right]  :=\left(  \frac{\left(  1-\rho\right)  ^{3}}{\rho^{2}%
-\rho+1},\frac{\left(  1-\rho\right)  ^{2}}{\rho^{2}+1}\right]  \label{int}%
\end{equation}
and that it is possible to give an explicit characterization of the set of
convex quadrilaterals with maximal ratio $M$. The fact that $M-m$ has a
maximum value of about $.034$ and is usually much smaller explains the near
constancy of $r\left(  \rho,ABCD\right)  $ as $ABCD$ varies. Here are the
graphs of $M$ and $m$.%
\begin{center}
\includegraphics[
height=3.1566in,
width=4.0517in
]%
{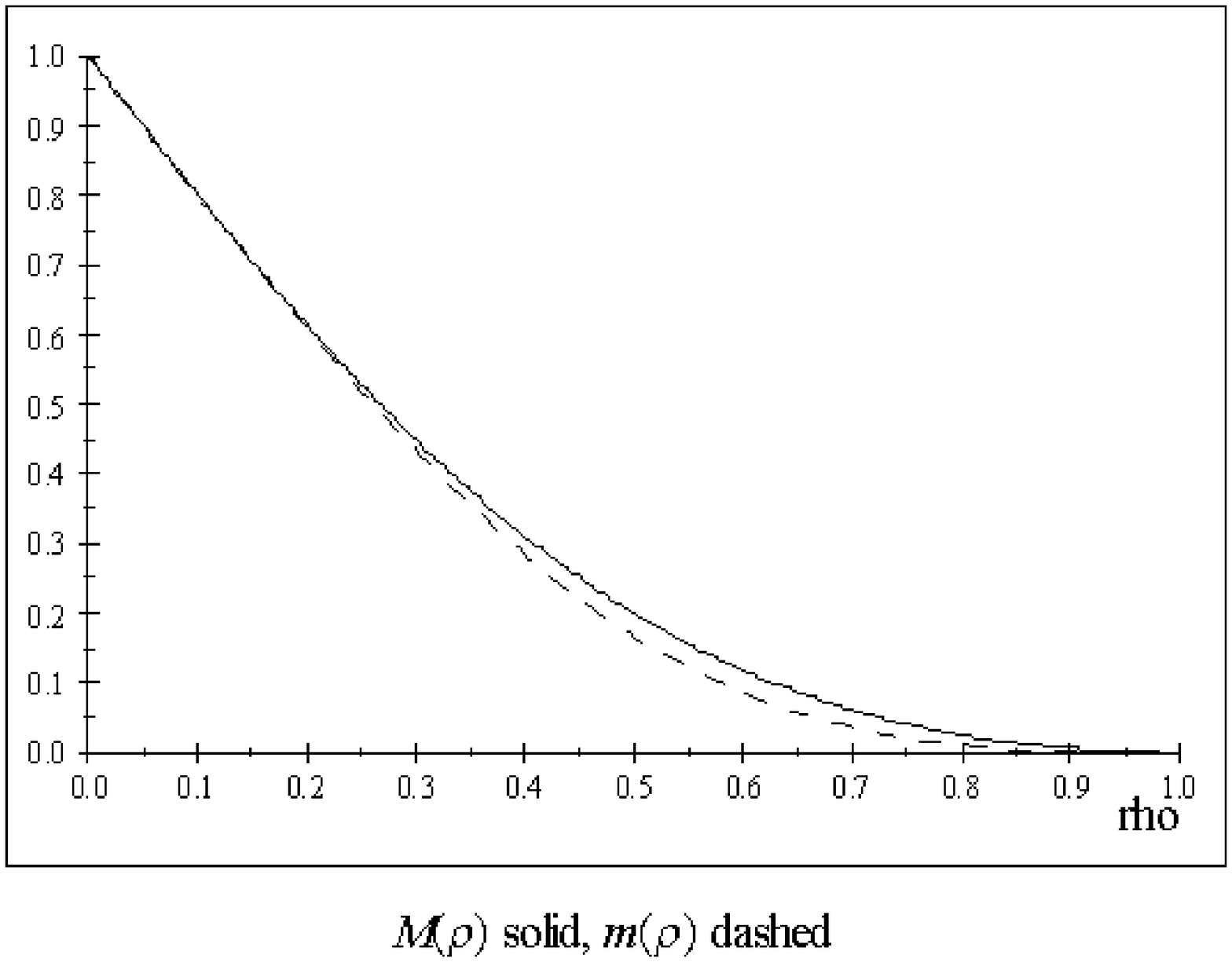}%
\end{center}
A more delicate look at the graph of $M-m=\allowbreak\frac{\rho^{3}\left(
\rho-1\right)  ^{2}}{\left(  \rho^{2}+1\right)  \left(  \rho^{2}%
-\rho+1\right)  }$ shows that as \textquotedblleft constant\textquotedblright%
\ as $r$ is in the original $\rho=1/2$ case, it is even \textquotedblleft more
constant\textquotedblright\ when $\rho$ is close to the endpoints $0$ and $1$.
(Actually the maximum value of $M-m$ of about $.034$ is achieved at the unique
real zero of $\rho^{5}-\rho^{4}+6\rho^{3}-6\rho^{2}+7\rho-3$ which is about
$.55$.)

The characterization proved in Theorem \ref{t:2} below shows that not only do
parallelograms have maximal ratio $M\left(  \rho\right)  $ for \textit{every}
$\rho$, but also they are the only quadrilaterals that have maximal ratio
$M\left(  \rho\right)  $ for more than one $\rho$.

\section{The midpoint case for parallelograms}

\begin{theorem}
\label{t:1}If each vertex of a parallelogram is joined to the midpoint of an
opposite side in clockwise order to form an inner quadrilateral, then the area
of the inner quadrilateral is one fifth the area of the original parallelogram.
\end{theorem}

\begin{proof}
In this picture,%
\begin{center}
\includegraphics[
height=2.124in,
width=5.4457in
]%
{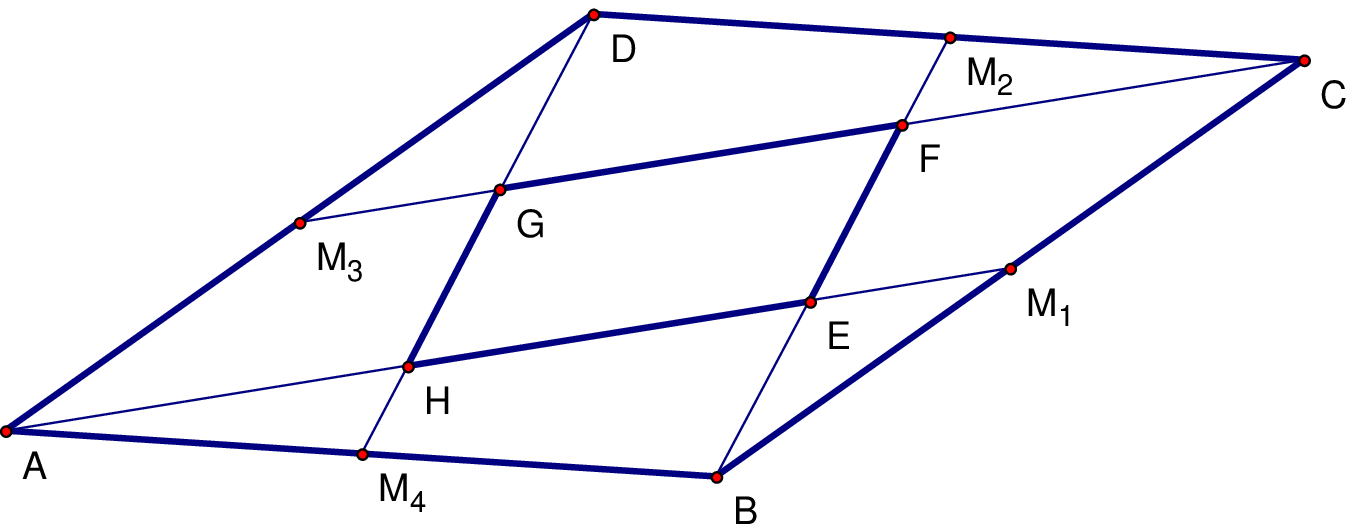}%
\end{center}
$ABCD$ is a parallelogram, and each $M_{i}$ is a midpoint of the line segment
it lies on. Cut apart the figure along all lines. Then by rotating clockwise
$180%
{{}^\circ}%
$ about point $M_{3}$, the reader can verify that we get $AHGG^{\prime}$(where
$G^{\prime}$ is the image of $G$ under the rotation) congruent to $EFGH$.
Similarly, each of the triangles $ABE$, $BCF$, and $CDG$ may be dissected and
rearranged to form a parallelogram, each congruent to $EFGH$. Thus, the pieces
of $ABCD$ can be rearranged into five congruent parallelograms, one of which
is $EFGH$, which therefore has area $1/5$ the area of $ABCD$.
\end{proof}

This result is a special case of Corollary \ref{c:3} below, but is included
because of the elegant and elementary nature of its proof.

\section{\label{s:3}The filling of $\left(  m,M\right]  $ and the
characterization}

\begin{theorem}
\label{t:2}Let $A,B,C,D$ be (counterclockwise) successive vertices of a convex
quadrilateral. Define $EFGH$ as the inner quadrilateral formed by joining
vertices to $\rho$-points as described in Section \ref{s:1}. Construct point
$P$ so that $ABCP$ is a parallelogram. Locate(as in the following figure)
$C^{\prime}$ and $C^{\prime\prime}$ on $\overrightarrow{BC}$ so that
$C^{\prime}$ is a distance $\rho BC$ from $C$ and $C^{\prime\prime}$ is a
distance $\left(  1/\rho\right)  BC$ from $C$ and let $S=\overrightarrow
{C^{\prime}P}\cup\overrightarrow{C^{\prime\prime}P}$. Then the ratio $r$
defined by equation (\ref{ratio}) is maximal exactly when $D$ is on $S^{\ast
}=S\cap\operatorname*{int}\left(  \angle ABC\right)  \cap\operatorname*{ext}%
\left(  \triangle ABC\right)  $. Furthermore the set of possible ratios is%
\[
\left(  \frac{\left(  1-\rho\right)  ^{3}}{\rho^{2}-\rho+1},\frac{\left(
1-\rho\right)  ^{2}}{\rho^{2}+1}\right]  .
\]

\end{theorem}

In the figure below, $S^{\ast}$ is indicated by the thickened portions of the
rays composing $S$.%
\begin{center}
\includegraphics[
trim=-1.533271in 0.000000in -1.533271in 0.000000in,
height=2.1223in,
width=5.1923in
]%
{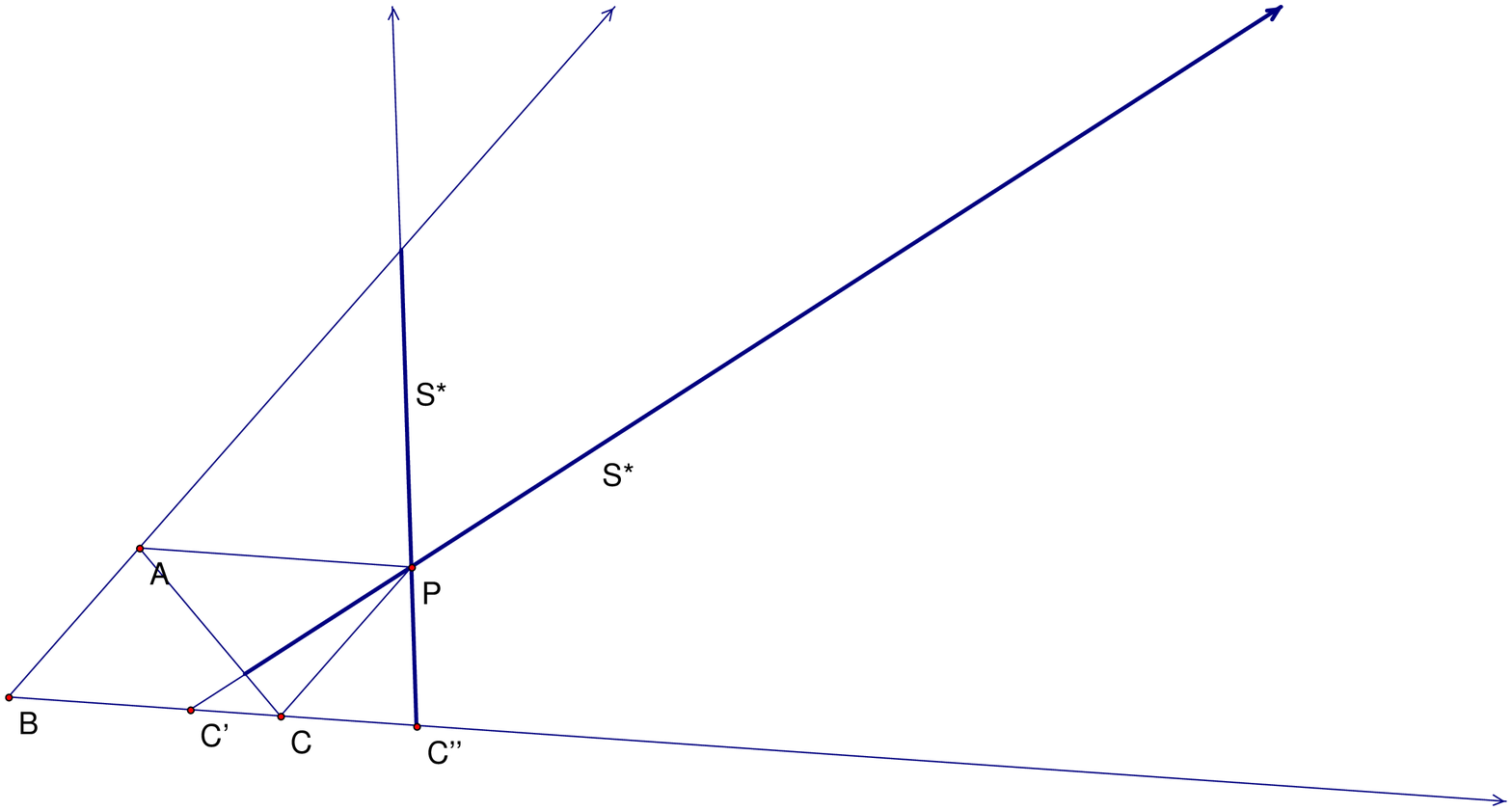}%
\end{center}

\begin{proof}
Fix $\rho$ and apply an transformation that maps $A,B,C,D$ successively to
$\left(  0,1\right)  ,\left(  0,0\right)  ,\left(  1,0\right)  ,\left(
x,y\right)  $. Since an affine transformation preserves both linear length
ratios and area ratios, it is enough to prove the theorem after the
transformation has been applied. Observe that $P$ has become $\left(
1,1\right)  $, and the image of $S$ has become a pair of perpendicular rays
through $\left(  1,1\right)  $ with slopes $\rho$ and $-1/\rho$. Here is the
situation.%
\begin{center}
\includegraphics[
height=3.6668in,
width=4.8983in
]%
{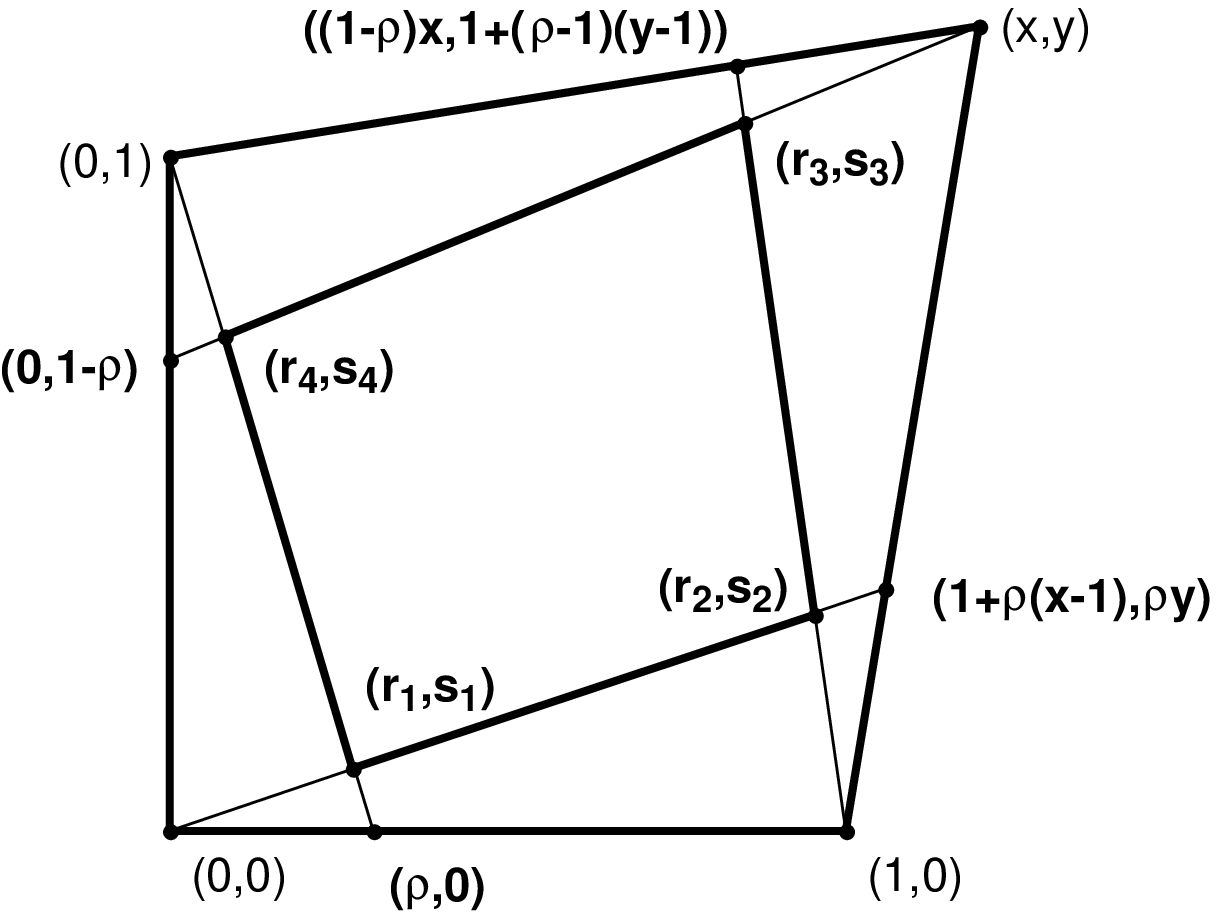}%
\end{center}
The line from $\left(  0,0\right)  $ to $\left(  x,y\right)  $ divides the
outer quadrilateral$\ $into two triangles, one of area $x/2$ and the other of
area $y/2$, so that its area is $\left(  x+y\right)  /2$. To find the area of
the inner quadrilateral, we first determine $\{r_{1},s_{1},\dots,s_{4}\}$ in
terms of $x$, $y$ and $\rho$ by equating slopes. For example, the equations
\[%
\begin{array}
[t]{l}%
\dfrac{s_{1}-0}{r_{1}-0}=\dfrac{\rho y-0}{1+\rho\left(  x-1\right)  -0}\\
\dfrac{s_{1}-0}{r_{1}-\rho}=\dfrac{1-0}{0-\rho}%
\end{array}
\]
can easily be solved for $r_{1}$ and $s_{1}$. The area of the interior
quadrilateral is%
\[
\frac{\left(  r_{1}s_{2}-r_{2}s_{1}\right)  +\left(  r_{2}s_{3}-r_{3}%
s_{2}\right)  +\left(  r_{3}s_{4}-r_{4}s_{3}\right)  +\left(  r_{5}s_{1}%
-r_{1}s_{5}\right)  }{2}.
\]
This is the $n=4$ case of a well-known formula for the area of an
$n$-gon\cite{Bra} which can be proved by first proving the formula for
triangles and then using induction, or by using Green's Theorem. Some computer
algebra produces this formidable and seemingly intractable formula for
$r\left(  x,y\right)  $.%
\[
\frac{\left(  \rho-1\right)  ^{2}\left(
\begin{array}
[c]{c}%
\rho^{4}y^{4}-\rho^{3}y^{4}-3\rho^{5}xy^{3}+2\rho^{4}xy^{3}+\rho^{3}%
xy^{3}-2\rho^{2}xy^{3}\\
+2\rho^{5}y^{3}-6\rho^{4}y^{3}+\rho^{3}y^{3}+2\rho^{2}y^{3}-\rho y^{3}%
+\rho^{6}x^{2}y^{2}\\
-\rho^{5}x^{2}y^{2}-6\rho^{4}x^{2}y^{2}+4\rho^{3}x^{2}y^{2}-\rho^{2}x^{2}%
y^{2}-\rho x^{2}y^{2}\\
-3\rho^{6}xy^{2}+10\rho^{5}xy^{2}+3\rho^{4}xy^{2}-13\rho^{3}xy^{2}+5\\
\ast\rho^{2}xy^{2}+\rho xy^{2}-xy^{2}+\rho^{6}y^{2}-7\rho^{5}y^{2}+6\rho
^{4}y^{2}+7\\
\ast\rho^{3}y^{2}-7\rho^{2}y^{2}+2\rho y^{2}+2\rho^{5}x^{3}y-2\rho^{4}%
x^{3}y-3\rho^{3}\\
\ast x^{3}y+2\rho^{2}x^{3}y-\rho x^{3}y-3\rho^{6}x^{2}y-5\rho^{5}x^{2}%
y+15\rho^{4}\\
\ast x^{2}y-\rho^{3}x^{2}y-7\rho^{2}x^{2}y+4\rho x^{2}y-x^{2}y+5\rho^{6}xy\\
-3\rho^{5}xy-21\rho^{4}xy+18\rho^{3}xy-\rho^{2}xy-3\rho xy+x\\
\ast y-2\rho^{6}y+5\rho^{5}y+4\rho^{4}y-12\rho^{3}y+6\rho^{2}y-\rho y\\
+\rho^{4}x^{4}-\rho^{3}x^{4}-3\rho^{5}x^{3}-2\rho^{4}x^{3}+5\rho^{3}%
x^{3}-2\rho^{2}x^{3}\\
+\rho^{6}x^{2}+8\rho^{5}x^{2}-6\rho^{4}x^{2}-5\rho^{3}x^{2}+5\rho^{2}%
x^{2}-\rho x^{2}\\
-2\rho^{6}x-5\rho^{5}x+12\rho^{4}x-4\rho^{3}x-2\rho^{2}x+\rho x+\rho^{6}\\
-4\rho^{4}+4\rho^{3}-\rho^{2}%
\end{array}
\right)  }{\left(
\begin{array}
[c]{c}%
\left(  y+\rho^{2}x-\rho x+x+\rho-1\right)  \left(  \rho y+x+\rho^{2}%
-\rho\right) \\
\ast\left(  \rho^{2}y+\rho x-\rho+1\right) \\
\ast\left(  \rho^{2}y-\rho y+y+\rho^{2}x-\rho^{2}+\rho\right)
\end{array}
\right)  }%
\]
\ Convexity means that $\left(  x,y\right)  $ is constrained to the open
\textquotedblleft northeast corner\textquotedblright\ of the first quadrant
bounded by $Y\cup T\cup X,$ $Y=\left\{  \left(  0,y\right)  :y\geq1\right\}
$, $T=\left\{  \left(  x,1-x\right)  :0\leq x\leq1\right\}  $, $X=\left\{
\left(  x,0\right)  :x\geq1\right\}  $. Restricting $r$ to $Y,$ we get a
formula for $r(y)=r\left(  0,y\right)  $. Taking the derivative unexpectedly
gives this simple, completely factored formula:%
\[
r^{\prime}\left(  y\right)  =\frac{\left(  \rho-1\right)  ^{2}\rho^{5}\left(
y-\frac{\rho+1}{\rho}\right)  \left(  y-\left(  \rho-1\right)  \right)
}{\left(  \rho^{2}y-\left(  \rho-1\right)  \right)  ^{2}\left(  \left(
\rho^{2}+1-\rho\right)  y-\rho\left(  \rho-1\right)  \right)  ^{2}}.
\]
The only solution to $r^{\prime}\left(  y\right)  =0$ with $y\in Y$ has
$y=\frac{\rho+1}{\rho}$. It quickly follows that on $Y$, $r$ attains a maximum
value of $M$ at $\left(  0,\frac{\rho+1}{\rho}\right)  $ and is minimized by
$m$ at the endpoints $\left(  0,1\right)  $ and $\left(  0,\infty\right)  $.
(By this we mean that $\lim_{y\rightarrow\infty}r\left(  0,y\right)  =m$.)
Similarly on $T,$ the derivative of $r\left(  x\right)  =r\left(
x,1-x\right)  $ has the following fully factored form%
\[
r^{\prime}\left(  x\right)  =\frac{\partial}{\partial x}r\left(  x,1-x\right)
=\frac{\left(  \rho+1\right)  \left(  \rho-1\right)  ^{2}\rho^{3}\left(
\left(  \rho-1\right)  x-\rho\right)  \left(  x-\frac{\rho}{\rho+1}\right)
}{\left(  \left(  \rho-1\right)  x-\rho^{2}\right)  ^{2}\left(  \rho\left(
\rho-1\right)  x-\left(  \rho^{2}+\left(  1-\rho\right)  \right)  \right)
^{2}},
\]
so that $r$ has minimum value $m$ at the endpoints $\left(  0,1\right)  $ and
$\left(  1,0\right)  $ and maximum value $M$ at $\left(  \frac{\rho}{\rho
+1},\frac{1}{\rho+1}\right)  $; while on $X$,
\[
r^{\prime}\left(  x\right)  =\frac{\partial}{\partial x}r\left(  x,0\right)
=\frac{\left(  \rho-1\right)  ^{2}\rho^{3}\left(  \rho^{2}+1\right)  \left(
x-\left(  \rho+1\right)  \right)  \left(  x-\left(  1-\rho\right)  \right)
}{\left(  x+\rho\left(  \rho-1\right)  \right)  ^{2}\left(  \left(  \rho
^{2}+1-\rho\right)  x+\rho-1\right)  ^{2}},
\]
so that $r$ has minimum value $m$ at the endpoints $\left(  1,0\right)  $ and
$\left(  \infty,0\right)  $ and maximum value $M$ at $\left(  \rho+1,0\right)
$. Motivated by these results we now sweep$\mathcal{\ }$the region of
permissible values of $\left(  x,y\right)  $ by line segments with
$y$-intercept $\eta$ and slope $-1/\rho$. From
\begin{align*}
&  r^{\prime}\left(  x\right)  =\frac{\partial}{\partial x}r\left(
x,-\frac{1}{\rho}x+\eta\right)  =\\
&  \frac{\left(
\begin{array}
[c]{c}%
\left(  \rho-1\right)  ^{2}\rho^{6}\left(  \rho^{2}+1\right)  ^{2}\left(
\eta-\frac{\rho+1}{\rho}\right)  ^{2}\left(  x-\dfrac{\rho^{2}+\eta\rho-\rho
}{\rho^{2}+1}\right) \\
\left(
\begin{array}
[c]{c}%
\left(  \eta\rho^{3}+\rho^{3}-3\rho^{2}-\eta\rho+3\rho-1\right)  x\\
-\left(  \rho^{3}+2\eta\rho^{2}+\rho-\eta\rho^{3}-\eta^{2}\rho^{2}-2\rho
^{2}-\eta\rho\right)
\end{array}
\right)
\end{array}
\right)  }{\left(
\begin{array}
[c]{c}%
\left(  \rho+\eta-1\right)  \left(  \eta\rho^{2}-\rho+1\right)  \left(
\left(  \rho^{3}-\rho^{2}+\rho-1\right)  x+\rho^{2}+\eta\rho-\rho\right)
^{2}\\
\left(  \left(  \rho^{3}-\rho^{2}+\rho-1\right)  x+\eta\rho^{3}-\rho^{3}%
-\eta\rho^{2}+\rho^{2}+\eta\rho\right)  ^{2}%
\end{array}
\right)  }%
\end{align*}
it is clear that $\eta=\frac{\rho+1}{\rho}$produces one arm of the image of
$S$. Finally for all other $\eta$, $r$ has mound-shaped behavior with minimum
values on the coordinate axes and reaches a maximum of $M$ where $y=-\frac
{1}{\rho}x+\eta$ intersects the other arm.
\end{proof}

Recall that we have defined $\rho$-points in terms of counterclockwise
orientation. Although Theorem \ref{t:2} is true for clockwise orientation, we
stress that the value of $r$ depends, in general, on the orientation. In fact,
clockwise and counterclockwise orientations always give different values of
$r$ unless $D$ lies on the diagonal $\overrightarrow{BP}$.

Setting $\rho=1/2$ in Theorem \ref{t:2} yields this corollary.

\begin{corollary}
\label{c:2}Let $A,B,C,D$ be (counterclockwise) successive vertices of a convex
quadrilateral. Define $EFGH$ as the inner quadrilateral formed by joining
vertices to midpoints as described in Section \ref{s:1}. Construct point $P$
so that $ABCP$ is a parallelogram. Locate $C^{\prime}$ and $C^{\prime\prime}$
on $\overrightarrow{BC}$ so that $C^{\prime}$ is a distance $\left(
1/2\right)  BC$ from $C$ and $C^{\prime\prime}$ is a distance $2BC$ from $C$
and let $S=\overrightarrow{C^{\prime}P}\cup\overrightarrow{C^{\prime\prime}P}%
$. Then the ratio $r$ defined by equation (\ref{ratio}) is maximal exactly
when $D$ is on $S^{\ast}=S\cap\operatorname*{int}\left(  \angle ABC\right)
\cap\operatorname*{ext}\left(  \triangle ABC\right)  $. Furthermore the set of
possible ratios is%
\[
\left(  \frac{1}{5},\frac{1}{6}\right]  .
\]

\end{corollary}

Another corollary of Theorem \ref{t:2} is the following generalization of
Theorem \ref{t:1} from midpoints to $\rho$-points.

\begin{corollary}
\label{c:3}If each vertex of a parallelogram is joined to the $\rho$-point of
an opposite side in counterclockwise order to form an inner quadrilateral,
then the area of the inner quadrilateral is $\frac{\left(  1-\rho\right)
^{2}}{\rho^{2}+1}$ times the area of the original parallelogram.
\end{corollary}

A nice geometry exercise is to prove this corollary avoiding the calculus part
of the proof of Theorem \ref{t:2}. Hint: Performing the affine transformation
we may assume that the original quadrilateral is the unit square. Use slope
considerations to see that the interior quadrilateral is actually a rectangle.
Use length considerations to see that it is a square of side length
$\sqrt{\frac{\left(  1-\rho\right)  ^{2}}{\rho^{2}+1}}$.

\end{document}